\documentclass[reqno]{amsart}

\usepackage{amsmath,amssymb,amscd,amsthm,accents} \usepackage{graphicx}
\DeclareGraphicsExtensions{.eps} \usepackage{mathrsfs}
\usepackage[mathcal]{eucal}
\usepackage{enumitem}
\usepackage{geometry}

\newtheorem{Thm}{Theorem}{\bfseries}{\itshape}
\newtheorem*{Thm*}{Theorem}{\bfseries}{\itshape}
\newtheorem{Cor}{Corollary}{\bfseries}{\itshape}
{\bfseries}{\itshape}
{\bfseries}{\itshape}
\newtheorem*{Lem*}{Lemma}{\bfseries}{\itshape}
{\bfseries}{\itshape}
{\bfseries}{\itshape}
{\bfseries}{\rmfamily}
{\scshape}{\rmfamily}
\newtheorem{Rem}[Cor]{Remark}{\scshape}{\rmfamily}
{\bfseries}{\itshape}

\let\tildeaccent=\~ \let\hataccent=\^
\renewcommand\~[1]{\widetilde{#1}}

\def\<{\left<} \def\>{\right>} \def\({\left(} \def\){\right)}

 \def\norm#1{\left\Vert #1
  \right\Vert}

\let\polishL=l \def\Zoladek.{\.Zol\c adek}

\def\etc.{\emph{etc}.}

\def\:{\colon}    \def\N{{\mathbb N}}

\let\PolishL=\L 
\def\L{{\mathbb L}}

  \def\F{{\varPhi}}
 
 \def\d{\,\mathrm d}

 \def\Lojas.{\PolishL ojasiewicz}

\def\cal{\mathcal}
\def\d{{\mathrm d}}

\begin{document}

\title{Effective Andr\'e-Oort for non-compact curves in Hilbert
 modular varieties}
  
\author{Gal Binyamini}
\address{Weizmann Institute of Science, Rehovot, Israel}
\email{gal.binyamini@weizmann.ac.il}

\author{David Masser}
\address{Departement Mathematik und Informatik, Universit\"at Basel, Spiegelgasse 1, 4051 Basel, Switzerland}
\email{David.Masser@unibas.ch}

\thanks{This research was supported by the ISRAEL SCIENCE FOUNDATION
  (grant No. 1167/17) and by funding received from the MINERVA
  Stiftung with the funds from the BMBF of the Federal Republic of
  Germany. This project has received funding from the European
  Research Council (ERC) under the European Union's Horizon 2020
  research and innovation programme (grant agreement No 802107)}

\subjclass[2010]{11G10, 11G15, 11G18, 11G50}

\date{\today}

\begin{abstract}
  In the proofs of most cases of the Andr\'e-Oort conjecture, there
  are two different steps whose effectivity is unclear: the use of
  generalizations of Brauer-Siegel and the use of Pila-Wilkie. Only
  the case of curves in ${\bf C}^2$ is currently known effectively (by
  other methods).

  We give an effective proof of Andr\'e-Oort for non-compact curves in
  every Hilbert modular surface and every Hilbert modular variety of
  odd genus (under a minor generic simplicity condition). In
  particular we show that in these cases the first step may be
  replaced by the endomorphism estimates of W\"ustholz and the second
  author together with the specialization method of Andr\'e via
  G-functions, and the second step may be effectivized using the
  Q-functions of Novikov, Yakovenko and the first author.
\end{abstract}

\maketitle
\date{\today}

\section{Introduction}

In this note we discuss the known results on effective Andr\'e-Oort,
all of which are restricted to the context of modular curves, and
announce a new effective result in the context of Hilbert modular
varieties. We provide a sketch of the key ideas, and will provide full
details in a paper under preparation.

\subsection{The Andr\'e-Oort conjecture and effectivity}

In the Andr\'e-Oort conjecture one considers a suitable ambient space
$\cal X$ - in this note we will mention only
${\cal X}={\bf C}^n~(n \geq 2)$ or ${\cal X}={\cal A}_g~(g \geq 2)$
the moduli space of principally polarized abelian varieties of
dimension $g$ - equipped with irreducible algebraic subvarieties that
are known as special subvarieties. Then one takes an algebraic
subvariety $V$ in $\cal X$, and one expects that if the special points
of $V$ are Zariski-dense in $V$, then $V$ is itself special.

For ${\cal X}={\bf C}^n$ the special points are the $(j_1,\ldots,j_n)$
whose coordinates are values of the elliptic modular function at
quadratic points of the standard upper half-plane. A special variety
is defined essentially by relations $F(x_r,x_s)=0$ where the $F$ are
modular polynomials. The conjecture was proved for $n=2$ by Andr\'e
himself \cite{A2} and a high point was reached with Pila's paper \cite{P2} for
all $n$ (and even more).

For ${\cal X}={\cal A}_g$ the special points are the abelian varieties
(not necessarily simple) with complex multiplication CM, in the sense
that their endomorphism algebras contain a commutative algebra of
dimension $2g$ over $\bf Q$. For special varieties the description is
not so elementary, but when $g=2,3$ it can also be reduced to
properties of the endomorphism algebras. The conjecture was proved for
$g=2$ by Pila and Tsimerman \cite{PT} and then for all $g$ by
Tsimerman \cite{T2}. 

In both situations ${\cal X}={\bf C}^n$ or
${\cal X}={\cal A}_g$ the proofs remain ineffective in the sense that
when $V$ is given, there is in general no method to find the set of
its special points.

For ${\cal X}={\bf C}^n$ some effective results are known in
particular cases. The case of a general curve $V$ has been treated
independently by K\"uhne \cite{K1} and Bilu, Zannier and the second
author \cite{BMZ}, and this remains the only case where the
Andr\'e-Oort conjecture is known effectively. Since then some progress
has been achieved for $V$ of higher dimension: see Bilu and K\"uhne
\cite{BK} for the linear case, and the first author \cite{B2} for an
extension to arbitrary degree hypersurfaces under a genericity
assumption on the leading homogeneous part. However no general result
for say a surface $V$ is known.

\subsection{Statement of the main result}

From now on $V$ will be a curve; thus if $V$ is not itself special
then it contains at most finitely many special points. For
${\cal X}={\cal A}_g$ there are no effective results up to now. But we
can prove the following. For simplicity we assume that $V$ is defined
over $\overline{\bf Q}$; as special points are also defined over
$\overline{\bf Q}$ there is no problem to extend to general fields
provided these are ``effectively given'' as for example in Fr\"ohlich
and Shepherdson \cite{FS}.

\begin{Thm}
  Let $V$ be a curve in ${\cal A}_g$, defined over $\overline{\bf Q}$,
  that is not special. Suppose
  \begin{itemize}
  \item[(a)] $g=2$ or $g>2$ is odd,
  \item[(b)] $V$ lies in a Hilbert modular variety,
  \item[(c)] $V$ is not compact,
  \item[(d)] a generic point of $V$ corresponds to a simple abelian variety.
  \end{itemize}
  Then the finitely many special points on $V$ can be found effectively.
\end{Thm}

We may identify $V$ with a single abelian variety $A$, defined over
the function field of a curve $B$, that is not isotrivial. Then $(b)$
is equivalent to $A$ having real multiplication RM, in the sense
that its endomorphism algebra contains a totally real field of
degree $g$ over $\bf Q$. And $(c)$ is
equivalent to $A$ having bad reduction somewhere. Also $(d)$ is
equivalent to $A$ being simple. For all but finitely many points $b$
of $B$, there is a specialized abelian variety $A_b$, and our claim is
then that the finitely many $b$ such that $A_b$ has CM can be found
effectively.

It seems plausible that condition $(d)$ can be dropped. Aside from
this technical issue, our main theorem applies to every non-compact
curve in every Hilbert modular surface or Hilbert modular variety of
odd genus, giving an analog for the effective result for curves in
${\bf C}^2$ by \cite{K1,BMZ} (where every curve is non-compact). In
particular, in Hilbert modular surfaces every non-special curve is
generically simple, and in these cases condition (d) is indeed
superfluous\footnote{We thank Jonathan Pila for pointing out this
  special case to us.}. To our knowledge no similar effective results
are known for any compact families in ${\cal A}_g$, and treating such
a problem would presumably require very different ideas, as all known
effective proofs involve some form of asymptotic around a cusp.

For small genus we may express the fact that $V$ is not special also
in a more elementary way. Thus for $g=2,3$ we get the following
consequences.

\begin{Cor}
  Let $A$ be a simple abelian surface, defined over
  ${\overline{\bf Q}}(B)$, that is not isotrivial. Suppose that $A$
  has real multiplication but not quaternion multiplication, and bad
  reduction somewhere. Then the finitely many $b$ in $B$ such that
  $A_b$ has complex multiplication can be found effectively.
\end{Cor}

An example with $B={\bf P}_1$ is the Jacobian of
\begin{equation}\label{1}
y^2=x^5+(t-2)x^4-(2t-3)x^3+tx^2-2x+1,
\end{equation}
with real multiplication by ${\bf Q}(\sqrt{5})$. This comes from
specializing a two-dimensional family due to Wilson \cite{W}. 
\begin{Cor}
  Let $A$ be a simple abelian threefold, defined over
  ${\overline{\bf Q}}(B)$, that is not isotrivial. Suppose that $A$
  has real multiplication and bad reduction somewhere. Then the
  finitely many $b$ in $B$ such that $A_b$ has complex multiplication
  can be found effectively.
\end{Cor}

An example again with $B={\bf P}_1$ is the Jacobian of
\begin{equation}\label{2}
y^2 = x^7-tx^6+34x^5+(8t+156)x^4+48x^3-16tx^2-32x + 64,
\end{equation}
with real multiplication by ${\bf Q}(\cos(2\pi/7))$. This comes from
specializing a two-dimensional family due to Mestre \cite{Me}. 

\begin{Rem}
  Of course such results are impossible for $g=1$, because all
  elliptic curves have RM and among them are infinitely many
  non-isomorphic (and even non-isogenous) curves with CM. For example,
  there are infinitely many $t$ such that
  \begin{equation}\label{3}
    y^2=x(x-1)(x-t)
  \end{equation}
  has CM.  
\end{Rem}

In fact our arguments are not restricted to real
multiplication. For example, the Jacobian of the curve
\begin{equation}\label{4}
y^2=x(x-1)(x-t)(x-t^2)(x-t^4)
\end{equation} 
has endomorphism ring $\bf Z$ (see \cite{Ma} p.296) and we can effectively
find all complex $t$ for which it is simple with complex
multiplication. Similarly for
\begin{equation}\label{5}
y^2=x(x-1)(x-t)(x-t^2)(x-t^4)(x-t^5)(x-t^8).
\end{equation}

\subsection{G-functions and comparison to Andr\'e's results}

Already Andr\'e \cite{A1} had observed the applicability of
G-functions to problems connected to special points in Shimura
varieties. In particular, the theorem in 
\cite{A1} (p.201) implies a result similar to our main theorem, with the
caveat that one counts CM points $s$ only with bounded residual degree
$[{\bf Q}(s):{\bf Q}]$. This was one of the motivations leading up to Andr\'e's
formulation of his general conjecture - see in particular \cite{A1} (pp. 215,216). 
Unlike the subsequent proof of the general conjecture,
Andr\'e's argument is effective. Our main contribution is to eliminate
the condition of bounded residual degree in our situation while preserving the
effectivity of the argument.

The recent papers of Daw and Orr \cite{DO1,DO2} also
use G-functions to study instances of the Zilber-Pink
conjecture. Though they do not consider effectivity, it seems
plausible that effective versions of their results may also be
possible to derive using arguments similar to those employed in the
present note.

\begin{Rem}
  In fact \cite{A1} assumes that $g>2$ is odd, and we are pleased to
  thank Yves Andr\'e for helping us to see that his work can be
  extended to $g = 2$.
\end{Rem}

We may note that G-functions were used for a different purpose by K\"uhne in \cite{K2}.
\subsection{Overview of the proof}

We use the strategy introduced by Pila-Zannier in \cite{PZ},
contrasting lower bounds on Galois orbits against upper bounds for
rational points in definable sets. As mentioned, the proofs in
\cite{PT} and \cite{T2} are ineffective, and this comes from two
different sources. First, a Brauer-Siegel bound due to Tsimerman
\cite{T1}, and second the Pila-Wilkie bound \cite{PW} for rational
points on definable sets.

We avoid Brauer-Siegel using a different version of the discriminant
which can be controlled by refinements of ``linear forms in abelian
logarithms'' due to W\"ustholz and the second author based on
transcendence theory. These introduce an extra height dependence which
can be dealt with using ``global relations for G-functions'' due to
Andr\'e, also based on transcendence theory. It is for this that we
need conditions $(a)$ and $(c)$ of our Theorem. And the combination of
$(b)$ and $(c)$ leads to Andr\'e's condition of multiplicative
reduction. For more details see section~\ref{sec:effective-bs}.

Several effective versions of the Pila-Wilkie theorem have been
established in the literature. However, none of these are general
enough to allow the type of counting needed in our context. More
specifically, such results have been either restricted to the Pfaffian
class of functions (see Jones and Thomas \cite{JT} for example) which
is not known to contain period maps; or restricted to counting in
compact subsets \cite{B1}. In our context it is necessary to count in
the entire (non-compact) fundamental domain.  We overcome the
compactness issue by employing the theory of Q-functions due to
Novikov, Yakovenko and the first author. This theory produces
effective estimates for the number of zeros, and other analytic
complexity measures, which crucially hold ``all the way to the
cusp''. For more details see section~\ref{sec:effective-pw}.

\subsection{Contents of this note}

In section~\ref{sec:effective-bs} we describe our method for achieving
Galois orbit lower bounds effectively by combining G-function methods
and endomorphism estimates. In section~\ref{sec:reduction-to-pw} we
describe how combining these bounds with Pila-Wilkie type bounds gives
the main theorem. In section~\ref{sec:effective-pw} we sketch how we
obtain such bounds effectively using Q-functions. 

For simplicity of the presentation we restrict attention to the
specific case of (\ref{1}) (essentially the same description applies
to (\ref{2}) above). Thus we assume that $t$ is such that the Jacobian
$A_t$ has complex multiplication CM. It is easy to see that $t$ is an
algebraic number. Let $d$ be its degree and $h$ its absolute
logarithmic height.

\section{Effectively avoiding Brauer-Siegel.}
\label{sec:effective-bs}

As $A_t$ has CM, the main result of Tsimerman \cite{T1} (p.1091) now
implies that
\begin{equation}\label{6}
|{\rm Disc}(Z({\rm End}(A_t)))| \leq cd^\kappa
\end{equation}
for a suitable discriminant of the centre of the endomorphism ring
over $\bf C$, where $c$ and $\kappa$ are absolute. The constant $c$ is
ineffective. Indeed this is so even in the analogue for (\ref{3}) in
dimension $g=1$, where the left-hand side corresponds to a
discriminant of an imaginary quadratic field and the right-hand side
corresponds to a class number. In fact (\ref{6}) holds for all $g$ with
suitable $c,\kappa$, as shown by Tsimerman \cite{T2} (p. 385).  

We will modify (\ref{6}) and factorize it through two effective inequalities. First we prove that
\begin{equation}\label{7}
0<{\rm Disc_r(End}(A_t)) \leq c_1(h+d)^{\kappa_1},
\end{equation} 
for a slightly different ``polarized discriminant'' and without taking
the centre. This result comes directly from certain ``endomorphism
estimates'' of W\"ustholz and the second author \cite{MW} (pp. 642,650)
obtained through transcendence techniques (with linear forms in
abelian logarithms). It holds much more generally for any $g$ and with
no assumptions on the endomorphism ring. Second we prove
that
\begin{equation}\label{8}
h\leq c_2d^{\kappa_2},
\end{equation}
which comes from a slight extension of a specialization result of
Andr\'e \cite[p.201]{A1}, also obtained through transcendence
techniques (with global relations for G-functions). The latter was
originally only for odd dimension $g \geq 3$ (the effectivity is still unknown for
$g=1$), and with a different sort of assumption about the endomorphism
ring; furthermore there is a new assumption involving multiplicative
reduction. This is automatically satisfied in the case of real
multiplication RM as soon as there is bad reduction somewhere. For our
(\ref{1}) the latter is also automatic because $B={\bf P}_1$ (and in
fact it holds at $t=\infty$, where there is multiplicative reduction
for (\ref{4}) and (\ref{5}) too). It is the use of transcendence
techniques that makes both (\ref{7}) and (\ref{8}) effective.

It is relatively easy to compare the discriminants in (\ref{6}) and
(\ref{7}), provided $A_t$ is isogenous to a power of a simple abelian
variety, also an automatic consequence of RM (which explains the extra
simplicity assumption for (\ref{4}) and (\ref{5}) above). Then
(\ref{7}) and (\ref{8}) give the required effective version of
(\ref{6}), assuming multiplicative reduction.

This (\ref{6}) is combined with an estimate for the Siegel matrix
$\tau_t$ of $A_t$ when normalized to lie in a fundamental domain of
the Siegel upper half-space $\cal H$. As $A_t$ has CM, it is
well-known that the entries are algebraic numbers, and so we can speak
of a (non-logarithmic absolute) height $H(\tau_t)$. For this Theorem
3.1 (p.206) of \cite{PT} implies that
\begin{equation}\label{9}
H(\tau_t) \leq c|{\rm Disc}(Z({\rm End}(A_t)))|^\lambda
\end{equation}
where $c,\lambda$ are absolute (here in dimension $g=2$) and
effective.

\section{Reduction to point-counting}
\label{sec:reduction-to-pw}

This (\ref{9}) is now set up for counting. The traditional way is to embed
$\cal H$ into ${\bf R}^6$ by taking the real and imaginary parts of
the three entries of the (symmetric) matrix. For our particular $t$
and $\tau_t$ it is easy to see that we get a point $P_t$ whose
coordinates are algebraic numbers of degree at most 16 with heights at
most $4H(\tau_t)^2$. If we allow $t$ for the moment to range over
$\bf C$ we get a real surface $Z$ in ${\bf R}^6$, actually
definable in the sense of \cite{PW}. Now by Theorem 1.6 of \cite{P1} (p.153) the number of our
particular $P_t$ with height at most $H$ which lie on the
transcendental part $Z^{\rm trans}$ is
\begin{equation}\label{10}N \leq c(\epsilon)H^\epsilon
\end{equation}
for any $\epsilon>0$. Here $c$ is up to now not effective. It is not
hard to show using Ax-Lindemann-Weierstrass as in Theorem 1.2 of \cite{PT}
(p.204) that $Z^{\rm trans}=Z$.

Now combining this with (\ref{6}) and (\ref{9}) we get
$N \leq c_1(\epsilon)d^{2\kappa\lambda\epsilon}$. A rather primitive
zero estimate (but also non-effective) yields a similar estimate for
the number of $t$. On the other hand we could have taken many
conjugates of just our single $t$, giving $N \geq d$. And now choosing
$\epsilon < 1/(2\kappa\lambda)$ shows that $d$ is bounded above. Thus
by (\ref{8}) so is $h$, and we finish with Northcott.  Thus we have to make
(10) effective.

\section{Effective Pila-Wilkie counting}
\label{sec:effective-pw}

The proof of (10) is based on the ideas of Bombieri and Pila \cite{BP}
and Pila and Wilkie \cite{PW} and involves two steps. For simplicity we will explain
the strategy for counting points only over $\bf Q$ - the same
ideas extend to points of any fixed degree, such as $16$ above. First, one
parametrizes the surface $Z$ using smooth charts with bounded
derivatives, and shows that in each of the charts the set of rational
points of height at most $H$ is interpolated by
$c(\epsilon) H^\epsilon$ algebraic hypersurfaces of degree
$C(\epsilon)$. Second, one consider the intersection of $Z$ with each
of these algebraic hypersurfaces separately and proceeds by induction.

It has been shown in \cite{B1} that this strategy can be
carried out effectively when the functions defining $Z$ satisfy
certain overdetermined systems of differential equations. In
particular, period maps fall within this class, and this result can
therefore be applied to our context. However, a crucial assumption in
 \cite{B1} is that one works in a compact domain where the
corresponding system of differential equations is non-singular. In our
context it is necessary to work with the non-compact set $Z$, and this
introduces significant technical difficulties, in particular in regard
to the parametrization step of the Pila-Wilkie proof.

\subsection{Q-functions}

To overcome the difficulty involving non-compact domains, we appeal to
the theory of \emph{Q-functions} developed by Novikov, Yakovenko and
the first author \cite{BNY1,BNY2}. Briefly, consider a system of
linear ordinary differential equations 
\begin{equation}\label{eq:conn-eq}
  \d X = \Omega\cdot X
\end{equation}
where $\Omega$ is a matrix of $\overline{\bf Q}$-rational meromorphic
one-forms on $V$. Suppose further that the singularities of the system
are regular, and that the monodromy is quasi-unipotent. Then, if $X$
denotes a fundamental solution for this system, each entry of $X$ is
called a Q-function. The \emph{complexity} of these Q-functions is
defined to be the maximum of the degrees and heights of the entries of
$\Omega$. By the direct-sum and tensor operations on connections,
Q-functions form an algebra and the complexity of sums and products is
suitably controlled in terms of the individual terms. In our
application, we will consider the equation~\eqref{eq:conn-eq}
associated to the Gauss-Manin connection for the $H^1$-cohomology of
the family of abelian varieties over $V$. Then $\tau_t$ is given by a
ratio of the two upper blocks of $X$ (thought of as a $2\times2$ block
matrix), and in particular each entry of $\tau_t$ is a ratio of
Q-functions (we call such ratios \emph{RQ-functions}). For example, the classical case 
(\ref{3}) corresponds to $\Omega=\begin{pmatrix}0&1 \\{1/4 \over t(1-t)}&{-1+2t\over t(1-t)}\end{pmatrix}{\rm d}t$ and $X=\begin{pmatrix}F&G \\F'&G'\end{pmatrix}$ for the classical hypergeometric function $F=F(1/2,1/2,1;t)$ and its derivative (see below for $G$).

Fix an \'etale coordinate $z:V\to{\bf C}$, and denote by
$\Sigma\subset{\bf C}$ its (finite) ramification locus. We may think
of the Q-functions as (multivalued) functions in the $z$-plane. The
main result of \cite{BNY1} shows that the number of zeros of a
Q-function in any disc $D\subset{\bf C}\setminus\Sigma$ is explicitly
bounded in terms of the complexity; and moreover, the same holds for a
punctured disc $D_\circ\subset{\bf C}\setminus\Sigma$ if one makes a
branch cut along some straight real line and counts zeros in a single
branch. This easily extends to RQ-functions. It is this ability to
count ``all the way to the cusp'' that we employ to avoid the
compactness restrictions.

\subsection{The point-counting strategy}

Another advantage of working with the complex analytic theory of
Q-functions in place of general real sets of o-minimal structures is
that we can view $Z$ as a one-dimensional \emph{complex} curve rather
than a real surface. Thus the Pila-Wilkie induction mentioned above
terminates after a single step, and the main issue is to suitably
parametrize $Z$, or at least the part of $Z$ that contains points up
to height $H$.

Our strategy is to cut up ${\bf C}\setminus\Sigma$ into finitely many
(punctured) discs and annuli such that in each piece, we can indeed
interpolate all points of height at most $H$ by $O(H^\epsilon)$ algebraic
curves of degree $C(\epsilon)$. Our basic technical tool is the
following theorem.

\begin{Thm}\label{thm:p-valent-interpolation}
  Let $f_1,f_2$ holomorphic on a disc in ${\bf C}$ and $p$-valent
  there. Then for $z$ in the concentric disc of half the
  radius, the points $(f_1(z),f_2(z))$ of height at most $H$ can be
  interpolated by $c(p,\epsilon)H^\epsilon$ curves of degree
  $C(\epsilon)$ for some explicit constants
  $c(p,\epsilon),C(\epsilon)$.
\end{Thm}

Theorem~\ref{thm:p-valent-interpolation} is similar to the main lemma
of Bombieri-Pila \cite{BP} (p.343), which has appeared in numerous
variations in the literature. However, the constants $C,c$ normally
depend on the maximum of the functions $f_1,f_2$ and their derivatives
-- and we replace this by the valency of these functions. Since
RQ-functions are known to be $p$-valent with $p$ effective, we can
apply this to any disc $D\subset{\bf C}\setminus\Sigma$ such that the
concentric disc of twice the radius does not meet $\Sigma$ without
having to control derivatives. In this way we cover the part of the
curve that lies away from the singularities $\Sigma$.

A different strategy is required for studying punctured discs
$D_\circ$ around the singular points. It is well-known that sections
of regular-singular connections with quasi-unipotent monodromy admit
converging expansions of the form
\begin{equation}\label{eq:asymp}
  f(z) = \sum_{\lambda\in\Lambda} z^\lambda P_\lambda(\log z)
\end{equation}
where $\Lambda\subset{\bf R}$ is a union of finitely many $\N$-cosets,
and each $P_\lambda$ is a polynomial of degree bounded by some fixed
integer. The expansion holds in a branch defined on $D_\circ$ with a
branch cut along a straight real line.

For example with (\ref{3}) these can be taken as essentially as $F$
and $G=iF\log z+\tilde F$ where $F=F(1/2,1/2,1;z)$ is as above and $\tilde F$ is a related power series in
$z$.

With the theory of Q-functions one has good control over the rate of
asymptotic convergence in~\eqref{eq:asymp}. Thus we can further
subdivide $D_\circ$ into discs and annuli such that
Theorem~\ref{thm:p-valent-interpolation} applies to all the discs; and
in each remaining annulus one of the terms
$z^\lambda P_\lambda(\log z)$ strongly dominates (the sum of) all
remaining terms.

The same applies for RQ-functions, with the dominant term now given in
the form $z^{\lambda} R_\lambda(\log z)$ for a rational function
$R_\lambda$. Suppose now that $f_1,f_2$ are two RQ-functions on an
annulus $E\subset D_\circ$, with dominant terms
$z^{\lambda_j}R_{\lambda_j}(\log z)$ for $j=1,2$. If $\lambda_1\neq0$
or $\lambda_2\neq0$ then the corresponding $f_1$ or $f_2$ can take
values of height $H$ only in a sub-annulus $E'\subset E$ of
logarithmic width $O(\log H)$. In this case it is a simple matter to
cover $E'$ by $O(\log H)$ discs and reduce to
Theorem~\ref{thm:p-valent-interpolation}.

The remaining case is $\lambda_1=\lambda_2=0$. It is a small miracle
that in this particular case, the leading terms of $f_1,f_2$ are
algebraically dependent (both being rational functions of
$\log z$). Since the remaining terms in the asymptotic expansion are
much smaller than these leading terms, it turns out that
$(f_1,f_2)(E)$ is extremely close to an algebraic curve -- and in
this case one is able to use a different, more naive argument to prove
the interpolation result. Briefly, we look for an upper-bound for an
interpolation determinant, as in the Bombieri-Pila method. But here
the interpolation determinant vanishes identically on the leading
terms (as they are algebraically dependent), and the remaining terms
are so small that a trivial majoration suffices for our purposes.

\subsection{Under the rug: complications around branch cuts}

Our approach is to parametrize $Z$ by means of discs and punctured
discs in $D,D_\circ\subset{\bf C}\setminus\Sigma$. Recall that in the
punctured discs we make branch cuts along straight real lines. This
means that we count rational values in the set $\tau(D_\circ)$ (of
Siegel matrices) for such a branch. In general the image
$\tau(D_\circ)$ will not be contained in a single fundamental domain:
the straight lines in the $z$-coordinate do not exactly correspond to
the semialgebraic walls of the fundamental domain in $\cal H$. If a
point $\tau(z)$ does not belong to the standard fundamental domain
then it is not possible, {\it a priori}, to estimate its height in
terms of the corresponding discriminant.

\begin{Rem}
  It is true that each $\tau(D_\circ)$ will meet finitely many
  fundamental domains, for instance by the definability of the
  universal covering map of the Siegel modular variety in an o-minimal
  structure \cite{PS}. However, computing this number effectively for
  a general curve $V$ does not seem straightforward.
\end{Rem}

We overcome this problem in a different manner. We show instead that a
suitable norm $\norm{\tau(z)}$ is effectively $O(|\log z|)$; for this
we consider $\tau$ as a map between the hyperbolic Riemann surface
$D_\circ$ and the hyperbolic Siegel space $\cal H$ and use the
Schwarz-Pick lemma. If $z$ corresponds to a CM point then $|\log z|$
is bounded by a polynomial in the corresponding discriminant. Thus
$\norm{\tau(z)}$ is similarly bounded, and this allows us to suitably
estimate its height by a polynomial in the discriminant even if
$\tau(z)$ does not belong to the standard fundamental domain.

\end{document}